\documentclass{elsart}
\usepackage{amssymb,amsmath}
\newcommand{\Nr}[1]{(\ref{#1})}
\journal{Applied Numerical Mathematics}
\volume{53}
\pubyear{2005}
\newcounter{mylastpage}
\issue{2--4}
\usepackage{hyperref}
\makeatletter
\def\ps@copyright{%
 \let\@oddhead\@empty
 \let\@evenhead\@empty
 \def\@oddfoot{\small\slshape\hskip-5em
   Published in \@journal\ \@volume\ (\the\@pubyear)\ no.\ \@issue, pp.\ \ESpagenumber{firstpage}--\ESpagenumber{mylastpage},
\href{http://dx.doi.org/10.1016/j.apnum.2004.08.023}{doi: 10.1016/j.apnum.2004.08.023}}%
 \let\@evenfoot\@oddfoot}
\makeatother
\begin{document}

\begin{frontmatter}

 \title{Convergence of Runge-Kutta Methods Applied to Linear Partial Differential-Algebraic Equations
 }
 \author[Darmstadt]{K. Debrabant\corauthref{cor1}}
 \ead{debrabant@mathematik.tu-darmstadt.de}
 \author[Halle]{K. Strehmel}
  \ead{strehmel@mathematik.uni-halle.de}
 \address[Darmstadt]{Darmstadt University of Technology, Department of Mathematics,
  Schlo{\ss}gartenstra{\ss}e 7, D-64283 Darmstadt, Germany}
\address[Halle]{Martin-Luther-Universit\"{a}t Halle-Wittenberg, FB Mathematik and Informatik,
 Institut f\"{u}r Numerische Mathematik, Postfach, D-06099 Halle (Saale), Germany}
 \corauth[cor1]{Corresponding author.}
\begin{abstract}
We apply Runge-Kutta methods to linear partial differential-algebraic equations of the form
 $A\:u_t(t,x) +
B(u_{xx}(t,x)+ru_x(t,x))+Cu(t,x) = f(t,x)$,
  where $A,B,C\in\Rset^{n,n}$ and the
matrix $A$ is singular.
We prove that under certain conditions the temporal convergence order of the
fully discrete scheme depends on the time index of the partial
differential-algebraic equation. In particular, fractional orders
of convergence in time are encountered. Furthermore we show that
the fully discrete scheme suffers an order reduction caused by the
boundary conditions. Numerical examples confirm the theoretical
results.
\end{abstract}

\end{frontmatter}

\section{Introduction}
In this paper we consider linear partial differential-algebraic equations (PDAEs) of the form
\begin{eqnarray}\label{pdaglgsys1}
 A\:u_t(t,x) +B\:\left( u_{xx}(t,x)+ru_{x}(t,x)\right) + C\: u(t,x) = f(t,x),
\end{eqnarray}
where $t\in (t_0 , t_e)$, $x \in \Omega=(-l,l)
\subset\Rset,$
 $A,B,C \in \Rset^{n,n}$ are constant matrices, $r\in\Rset,$
 $u,f: [t_0 , t_e]\times \overline{\Omega}\to \Rset^{n}$. We are interested in cases where the matrix $A$ is singular.
 The singularity of $A$ leads to the differential-algebraic aspect.\\
 It will always be tacitly assumed that the exact solution is as often differentiable as the numerical
 analysis requires.

 In contrast to parabolic initial boundary value problems with regular matrices $A$ and $B$, here we cannot
  prescribe initial and boundary values for all components of the solution vector,
 they have to fulfill certain consistency
conditions. We consider one example:
\begin{exmp}Superconducting coil (see Marszalek/Trzaska, Campbell/Marszalek \cite{MarszalekTrzaska,CampbellMarszalek.b}):
\[
\begin{pmatrix}{~\;0}&0\\{-\frac{LC}{l^2}}&{\frac{L}D}
\end{pmatrix}
u_{tt}-u_{xx}
+\begin{pmatrix}0&1\\0&0
\end{pmatrix}
u=0
\]
with $~x\in(0,l),~t>0$. $u_1(t,x)$ denotes the voltage, $u_2(t,x)$ denotes the divergence of the electric
field strength within the coil.
 $l$ is the length of the whole winding. $L$, $C$ and $D$ are further coil parameters.
Transformation to a partial differential-algebraic system of first order in $t$ yields
\begin{equation}\label{GlgSystemSupraMagnetSpule}
\begin{pmatrix}
0 & 0 &0 &0\\ 0 & 0 &-\frac{LC}{l^2}~~ & \frac LD\\ {1} & {0} & {0} & {0}\\ {0} & {1} & {0} & {0}
\end{pmatrix}
u_t-
\begin{pmatrix}
1 & 0 &0 &0\\ 0 & 1 & 0 &0\\ {0} & {0} & {0} & {0}\\ {0} & {0} & {0} & {0}
\end{pmatrix}
u_{xx} +
\begin{pmatrix}
0 & 1 &~\;0 &~\;0\\ 0 & 0 & ~\;0 &~\;0\\ {0} & {0} & {-1} & {~\;0}\\ {0} & {0} & {~\;0} & {-1}
\end{pmatrix}
u=0.
\end{equation}
As initial conditions we choose
\[u_1(0,x)=\left(\frac{E}l-\frac{CDEl}6\right) x+\frac{CDE}{6l}x^3,~u_3(0,x)=0\]
and as boundary conditions
\[u_1(t,0)=u_2(t,0)=0,~u_1(t,l)=E,~u_2(t,l)=CDE,\]
 where $E$ is the energizing
source voltage at the input of the coil.\\
 As the boundary values of $u_1$ and $u_2$ are constant,
 we get from the third and fourth equation of \Nr{GlgSystemSupraMagnetSpule} that $u_3$ and $u_4$ fulfill
 homogeneous boundary conditions.
  From the initial condition of $u_1$ and the first equation we derive
   \mbox{$u_2(0,x)=\frac{CDE}{l}x$.} With $u_3(0,x)=0$ and the third equation it follows
   \mbox{$u_{1t}(0,x)=0$} and therefore \mbox{$u_{1xxt}(0,x)=0$}. With the first equation this implies
   \mbox{$u_{2t}(0,x)=0$,}
    and with the fourth equation we get finally $u_{4}(0,x)=0$.\\
  Here we have chosen the prescribed initial and boundary values such that
  all initial and boundary values are compatible.
\end{exmp}
For further examples considering the determination of the initial and boundary values which cannot be prescribed see
Lucht/S./Eichler-Liebenow \cite{LuStreEiL}.

In the following we assume that for the numerical computation all initial values
\[
u(t_0,x)= \varphi(x),\;x\in\bar{\Omega},
\]
and all boundary values entering into the space discretization are known,
\begin{alignat*}4
B u(t,x)&=\;&&\psi(t,x),~&&x\in\partial\Omega, ~&&t\in [t_0 , t_e],
\end{alignat*}
where we restrict ourselves to Dirichlet boundary conditions to simplify the presentation.\\
Investigations of the convergence of Runge-Kutta methods applied to abstract parabolic differential equations
 can be found for example in
 Brenner/Crouzeix/Thom\'{e}e \cite{BrennerCrouzeixThomee}, Lubich/Ostermann \cite{LubichOstermann1993} and
 Ostermann/Thalhammer \cite{OstermannThalhammer}. The approach used there cannot be
  carried forward directly to the class of problems considered here because
   the matrix $A$
is singular.

This paper is organized as follows: In Section \ref{Abschnitt2} we derive a semi-discrete system based on finite
  differences. The result is a
 method-of-lines-DAE (MOL-DAE).\\
  Section \ref{Abschnitt3} is devoted to the Runge-Kutta approximation of
  the MOL-DAE. Under a regular transformation, the MOL-DAE of dimension $nN$ is decoupled into $N$ systems
  of dimension $n$, where $N$ denotes the number of grid points on the $x$-axis.  Furthermore, a Weierstrass-Kronecker
  transformation is used to decouple each of these systems into an ODE-system and an algebraic system.
  We introduce the differential time index of the linear PDAE and give the Runge-Kutta approximation to these subsystems.\\
In Section \ref{Abschnitt4} we prove the convergence of $L$-stable Runge-Kutta discretizations with constant step sizes.
The attained order of convergence in time depends on the
    differential time index of the PDAE and on the boundary conditions (homogeneous or inhomogeneous)
     which enter into the space discretization.\\
Numerical experiments are finally presented in Section \ref{Abschnitt5}. We
illustrate our convergence results for the backward Euler method and the 3-stage Radau IIA method.
\section{Space discretization}\label{Abschnitt2}
The discretization in space of problem (\ref{pdaglgsys1}) by means of
finite-differences results in a differential-algebraic equation
(MOL-DAE)
\begin{equation}\label{ZMOL}
M\dot{U}=DU(t)+\tilde{F}(t), \quad t_0\leq t\leq t_e,
\end{equation}
 where $U(t)$ is an $Nn$-dimensional real vector
consisting of approximations to $u$ at the grid points. Here $N$
denotes the number of grid points on the $x$-axis. The matrix $M$
is given by $M=I_N\otimes A$ and the matrix $D$ originates from
the discretization of the differential operator
$B\frac{\partial^2}{\partial x^2}$ by second order
difference-approximations, from the discretization of the
differential operator $rB\frac{\partial }{\partial x}$ by second
($\delta=\frac12$) or first order difference-approximations
($\delta\in[0,1]\setminus\{\frac12\}$) and from the matrix $C$,
i.e., $D$ is given by
\[D=-\frac1{h^2}P\otimes{B}-I_N\otimes{C},
  \]
where $I_{N}$ is the $N$-dimensional identity matrix,
\[ P=\begin{pmatrix}
-(2-hr(1-2\delta))
&  1+hr\delta &        \\
1+hr(\delta -1)& -(2-hr(1-2\delta))& 1+hr\delta \\
    \dots    &  \dots  & \dots        \\
     \dots   & \dots   & 1+hr(\delta -1)&
-(2-hr(1-\delta))
\end{pmatrix}
\]
and $h=\frac{2l}{N+1}$ denotes the constant grid size.
 The $Nn$-dimensional real vector $\tilde{F}(t)$ arises from the
right hand side $f$ of (\ref{pdaglgsys1}) and the boundary values which enter into the discretization.\\ We denote by
$U_h(t)$ the restriction of $u(t,x)$ to the spatial grid and by $\alpha_h(t)$ the space truncation error defined by
\begin{equation}\label{alpha}
\alpha_{h}(t):=M\dot{U}_{h}(t) -D U_{h}(t)-\tilde{F}(t).
\end{equation}
By Taylor expansion of the exact solution we get
\begin{alignat}2\label{alpha2}
\alpha_h(t)&= h^{p_x} (I_{N}\otimes B)
\gamma_{h}(t) &\quad \text{with}\quad
\|\gamma_{h}(t)\|_\infty&\leq K,
\end{alignat}
 where $p_x\in\{1,2\}$ is the order of approximation of the space
 discretization and $K$ is a positive constant, i.e.,
\[\alpha_h(t)=\mathcal{O}(h^{p_x})\quad
\text{as}\quad h\to 0.\]
 Furthermore, we can show  that there exists a regular matrix $Q$ with
\begin{alignat}2\label{ShiGlg}
Q^{-1}\frac{1}{h^2}PQ=\text{diag}\{\lambda_{1},\dots,\lambda_{N}\},&\quad
\end{alignat}
where
\[
\lambda_{j}=-\frac{ 2-hr(1-2\delta)}{h^2}+2\frac{1+hr\delta  }{h^2}
 \sqrt{\frac{1+hr(\delta -1)}{1+hr\delta}}\cos\frac{j\pi}{N+1}.
\]
In the discrete $L_2$-norm we have
 \begin{equation}\label{CKonst}
\max\{\|Q\|,\|Q^{-1}\|\}\leq C_{1}
 \end{equation}
  with
 a positive constant $C_1$ independent of $h$. Therefore, in the following this norm is used.\\
\section{Runge-Kutta approximations}\label{Abschnitt3}
In order to numerically advance in time the solution of the
MOL-DAE (\ref{ZMOL}), we employ an $s$-stage Runge-Kutta method
\begin{alignat*}{2}
U_{m+1}^{(i)}&=U_m+\tau\sum_{i=1}^sa_{ij}K_{m+1}^{(i)},&~
MK_{m+1}^{(i)}&=DU_{m+1}^{(i)}+\tilde{F}(t_m+c_i\tau),~ 1\leq i\leq s,\\
U_{m+1}&=U_m+\tau\sum_{i=1}^sb_iK_{m+1}^{(i)},
\end{alignat*}
where $a_{ij},b_i,c_i\in\mathbb{R}$ are the coefficients of the
method and $\tau=\frac{t_e-t_0}{M_e}$ the time step size.\\ For
the investigation of the convergence of the method, it is useful
to introduce the Runge-Kutta matrix
$\mbox{$\mathfrak{A}$}=(a_{ij})_{ij=1}^s$ and  the vector notation
 ${1\kern-0.25em{\rm l}}_s=(1,\dots,1)^{\top}\in
\mathbb{R}^s,$  \mbox{$b=(b_1,\cdots,b_s)^{\top}.$}\\ Then, with
the Kronecker product, we obtain the compact scheme
\begin{subequations}\label{NLRKGlg}
\begin{alignat}1
U_{m+1}&=U_m+\tau\left( b^{\top}\otimes I_{Nn}\right) K_{m+1}\label{NLRKGlg1},\\
S_{m+1}&={1\kern-0.25em{\rm l}}_s\otimes U_m+\tau\left(\mbox{$\mathfrak{A}$}\otimes I_{Nn}\right) K_{m+1}\label{NLRKGlg2},\\
\left( I_s\otimes M\right) K_{m+1}&=\left( I_s\otimes D\right) S_{m+1}+\bar{F}(t_{m+1}),\label{NLRKGlg3}~m=0,\dots,M_e-1,
\end{alignat}
\end{subequations}
where
$S_{m+1}=\left({U_{m+1}^{(1)}}^{\top},\dots,{U_{m+1}^{(s)}}^{\top}\right)^{\top},$
$K_{m+1}=\left({K_{m+1}^{(1)}}^{\top},\dots,{K_{m+1}^{(s)}}^{\top}\right)^{\top}$
and
$\bar{F}(t_{m+1})=\left({\tilde{F}\left(
t_m+c_1\tau\right)}^{\top},\dots,{\tilde{F}\left(
t_m+c_s\tau\right)}^{\top}\right)^{\top}.$

By the regular transformation \Nr{ShiGlg}, the MOL-DAE \Nr{ZMOL} can be decoupled into $N$
DAEs
\begin{equation}\label{DASys_entk}
A\dot{U}_{Qk}(t)=D_kU_{Qk}(t)+\tilde{F}_{Qk}(t),\quad k=1,\dots,N,
\end{equation}
with $D_k=-{\lambda_k} B-C$ and
\begin{alignat*}{1}
\left( U_{Q1}(t)^\top,\dots,U_{QN}(t)^\top\right)^\top&=\left( Q^{-1}\otimes I_n\right) U(t),\quad\\
\left( \tilde{F}_{Q1}(t)^\top,\dots,\tilde{F}_{QN}(t)^\top\right)^\top&=\left( Q^{-1}\otimes I_n\right)\tilde{F}(t).
\end{alignat*}
In the following we assume that the matrix pencil $\{D+\lambda M\}$, ${\lambda\in\Cset}$, is regular, which is equivalent to the
regularity of all the matrix pencils $\{ D_k+\lambda A\}$.
\begin{defn}
Suppose that all matrix pencils $\{D_k+\lambda A\}$, $k=1,\dots,N$, are regular and have the same index $\nu_{dt}$. Then
the differential time index of the linear PDAE \Nr{pdaglgsys1} is defined to be $\nu_{dt}$.
\end{defn}
According to Weierstrass and Kronecker there exist
 regular matrices $P_{k}$ and $Q_{k}$ with
\begin{subequations}\label{WKTransf}
\begin{eqnarray}
P_{k}AQ_{k}&=&\mathrm{diag}\{I_{n_{k 1}},\dots,I_{n_{k s_k}},N_{m_{k 1}},\dots,N_{m_{k l_k}}\},\\
P_{k}D_{k}Q_{k}&=&\mathrm{diag}\{R_{k 1},\dots,R_{k s_k},I_{m_{k 1}},\dots,I_{m_{k l_k}}\}\label{Dkvekdiag},
\end{eqnarray}
 where
\begin{equation}\label{RkiNkiDef}
R_{k i}=
  \left(\begin{array}{cccc}
      \varkappa_{k i} &  1 &    &   0    \\
        & \ddots & \ddots      &       \\
        &    & \varkappa_{k i} &    1   \\
    0    &    &       & \varkappa_{k i}
   \end{array}
\right)\in\Cset^{n_{k i},n_{k i}},~
N_{m_{k i}}=
  \left(\begin{array}{cccc}
      0 &  1 &    &     0    \\
        & \ddots & \ddots    &     \\
        &    & 0 &    1   \\
      0 &    &       & 0
      \end{array}
\right)\in\Cset^{m_{k i},m_{k i}}
 \end{equation}
\end{subequations}
(see Hairer/Wanner \cite{HairerWanner}),
 and for the differential time index of the PDAE it follows
\[\nu_{dt}=\max\limits_{k}\{m_{k i}:i=1,\dots,l_k\}.\]
Therefore, DAE \Nr{DASys_entk} is decoupled into systems of the form
\begin{subequations}
\begin{alignat}3\label{WK1}
\dot{U}_{1kl}(t)&=R_{kl}U_{1kl}(t)+\Tilde{F}_{1kl}(t),\quad &l&=1,\dots,s_k,\\
N_{m_{kl}}\dot{U}_{2kl}(t)&=U_{2kl}(t)+\Tilde{F}_{2kl}(t),\quad &l&=1,\dots,l_k\label{WK2}
\end{alignat}
\end{subequations}
with
\[
\left( U_{1k1}(t)^\top,\dots,U_{1ks_k}(t)^\top,U_{2k1}(t)^\top,\dots,U_{2kl_k}(t)^\top\right)^\top
=Q_k^{-1}U_{Qk}
\]
and
\[
\left(\Tilde{F}_{1k1}(t)^\top,\dots,\Tilde{F}_{1ks_k}(t)^\top,\Tilde{F}_{2k1}(t)^\top,\dots,\Tilde{F}_{2kl_k}(t)^\top\right)^\top
=P_k\Tilde{F}_{Qk}.
\]
Similarly, DAE \Nr{alpha} can be transformed to
\begin{subequations}
\begin{alignat}3
\label{alphaWK1}
\dot{U}_{h1kl}(t)&=R_{kl}U_{h1kl}(t)+\Tilde{F}_{1kl}(t)+\alpha_{h1kl}(t),\quad &l&=1,\dots,s_k,\\
N_{m_{kl}}\dot{U}_{h2kl}(t)&=U_{h2kl}(t)+\Tilde{F}_{2kl}(t)+\alpha_{h2kl}(t),\quad
&l&=1,\dots,l_k.
\end{alignat}
\end{subequations}
Runge-Kutta methods are invariant under the transformations
\Nr{ShiGlg} and \Nr{WKTransf}. Therefore, to analyze convergence
 it is sufficient to apply them to
 systems of the form \Nr{WK1} and \Nr{WK2}.
  Application to \Nr{WK1} yields
\begin{subequations}\label{NLRKGlgWK1}
\begin{alignat}1
U_{1kl,m+1}&=U_{1kl,m}+\tau\left( b^{\top}\otimes I_{n_{kl}}\right) K_{1kl,m+1}\label{NLRKGlg1WK1},\\
S_{1kl,m+1}&={1\kern-0.25em{\rm l}}_s\otimes U_{1kl,m}+\tau\left(\mbox{$\mathfrak{A}$}\otimes I_{n_{kl}}\right) K_{1kl,m+1}\label{NLRKGlg2WK1},\\
K_{1kl,m+1}&=\left( I_s\otimes R_{kl}\right) S_{1kl,m+1}+\bar{F}_{1kl}(t_{m+1}),\label{NLRKGlg3WK1}~m=0,\dots,M_e-1,
\end{alignat}
\end{subequations}
and to \Nr{WK2}
\begin{subequations}\label{NLRKGlgWK2}
\begin{alignat}1
U_{2kl,m+1}&=U_{2kl,m}+\tau\left( b^{\top}\otimes I_{n_{kl}}\right) K_{2kl,m+1}\label{NLRKGlg1WK2},\\
S_{2kl,m+1}&={1\kern-0.25em{\rm l}}_s\otimes U_{2kl,m}+\tau\left(\mbox{$\mathfrak{A}$}\otimes I_{n_{kl}}\right) K_{2kl,m+1}\label{NLRKGlg2WK2},\\
\left( I_s\otimes N_{m_{kl}}\right) K_{2kl,m+1}&=S_{2kl,m+1}+\bar{F}_{2kl}(t_{m+1}),\label{NLRKGlg3WK2}~m=0,\dots,M_e-1.
\end{alignat}
\end{subequations}
Now we  start our convergence investigations.
\section{Convergence estimates}\label{Abschnitt4}
At first we introduce the  global (space-time discretization) error
$e_{m+1}$ and the residual (space-time discretization) errors
$\delta_{m+1}$ and $\Delta_{m+1}$ at the time level $t=t_{m+1}$.
\begin{defn} The global error $e_{m+1}$ at $t_{m+1}$ is defined by
\begin{equation*}
e_{m+1}:=U_h(t_{m+1})-U_{m+1}
\end{equation*}
and the residual errors $\delta_{m+1},\Delta_{m+1}$ are given by
\begin{subequations}
\begin{equation}
{\delta_{m+1}}:=U_{h}(t_m+\tau)-U_{h}(t_m)- \tau\left(
b^{\top}\otimes I_{Nn}\right)\hat{K}_{m+1},
 \label{ResGlg2}
\end{equation}
\begin{equation}
\label{ResGlg1} \Delta_{m+1} := \hat{S}_{m+1}-{1\kern-0.25em{\rm l}}_s\otimes
U_h(t_m)-\tau\left(\mbox{$\mathfrak{A}$}\otimes I_{Nn}\right)\hat{K}_{m+1},
\end{equation}
\end{subequations}
where $\hat{S}_{m+1}$ and $\hat{K}_{m+1}$ are defined by the exact
solution $U_h(t)$ of the PDAE, i.e.,
\begin{alignat*}1
\hat{S}_{m+1}&:=\left( U_h(t_m+c_1\tau)^{\top},\dots U_h(t_m+c_s\tau)^{\top}\right)^{\top},\\
\hat{K}_{m+1}&:=\left(\dot{U}_h(t_m+c_1\tau)^{\top},\dots\dot{U}_h(t_m+c_s\tau)^{\top}\right)^{\top}.
\end{alignat*}
\end{defn}
\begin{defn}
The discretization scheme (\ref{NLRKGlg}) is convergent of order
$(p_x,p^\star)$, if the global error satisfies
 \[\|  e_{m+1}\| =\mathcal{O}(h^{p_x})+\mathcal{O}(\tau^{p^\star})
\quad \mathrm{for} \quad (m+1)\tau=\mathrm{const.},~\tau,\, h\to 0,\]
whenever $u(t,x)$ is sufficiently often differentiable.
\end{defn}
With the components
 $e_{Qk,m+1}$ defined by
 \[\left( e_{Q1,m+1}^\top,\dots,e_{QN,m+1}^\top\right)^\top:=(Q^{-1}\otimes I_n)e_{m+1}\] and \Nr{CKonst}
 we obtain
 the estimate
 \begin{equation}\label{eqzsglg}
\frac1{C_1}\sqrt{h\sum_{k=1}^N\|e_{Qk,m+1}\|^2}\leq \|e_{m+1}\|\leq C_1\sqrt{h\sum_{k=1}^N\|e_{Qk,m+1}\|^2}.
 \end{equation}
Letting
$\alpha_{h 1kl, m+1}=\left(\alpha_{h 1kl} (t_m+c_1\tau)^{\top},\dots,\alpha_{h 1kl}(t_m+c_s\tau)^{\top}\right)^{\top}
$
 we get with \Nr{alphaWK1} and \Nr{NLRKGlg3WK1}
 \[
\hat{K}_{1kl,m+1}-K_{1kl,m+1}=
\left( I_s\otimes R_{kl}\right)\left(\hat{S}_{1kl,m+1}-S_{1kl,m+1}\right)+\alpha_{h1kl, m+1}.
\]
Using (\ref{NLRKGlg2WK1}), the transformed components
$e_{1kl,m}=U_{h1kl}(t_{m})-U_{1kl,m}$ of the global discretization error and  the transformed components
$\Delta_{1kl,m+1}$ of
(\ref{ResGlg1}) we obtain
\begin{equation*}
\hat{S}_{1kl,m+1}-S_{1kl,m+1}={1\kern-0.25em{\rm l}}_s\otimes e_{1kl,m}
+\tau\left(\mbox{$\mathfrak{A}$}\otimes I_{n_{kl}}\right)\left(\hat{K}_{1kl,m+1}-K_{1kl,m+1}\right)
+\Delta_{1kl,m+1}\label{StGlgDiff2}.
\end{equation*}
Combining the last two equations leads to
\begin{equation}\label{Stern1}
G(\tau R_{kl})\left(\hat{K}_{1kl,m+1}-K_{1kl,m+1}\right)= \left({1\kern-0.25em{\rm l}}_s\otimes
R_{kl}\right) e_{1kl,m}
+\left( I_s\otimes R_{kl}\right)\Delta_{1kl,m+1}+\alpha_{h1kl, m+1},
\end{equation}
where $G(z)=I_s-z\mbox{$\mathfrak{A}$}$.
\begin{rem}
For a (matrix-valued) function $f(z):\Cset\to\Cset^{m,n}$ which is analytic in a neighbourhood of $\kappa_{ki}$, the matrix function
$f(R_{ki})$ with $R_{ki}$ given in \Nr{RkiNkiDef} is defined by
 (see Golub/van Loan \cite{GolubvanLoan})
\[
f(R_{ki}):=
  \left(\begin{array}{cccc}
      f_{i_1i_2}(\varkappa_{k i}) &  \frac{f_{i_1i_2}^{(1)}(\varkappa_{k i})}{1!} &  \dots  &   \frac{f_{i_1i_2}^{(n_{ki}-1)}(\varkappa_{k i})}{(n_{ki}-1)!}    \\
        & \ddots & \ddots      & \vdots      \\
        &    & f_{i_1i_2}(\varkappa_{k i}) &    \frac{f_{i_1i_2}^{(1)}(\varkappa_{k i})}{1!}   \\
    0    &    &       & f_{i_1i_2}(\varkappa_{k i})
   \end{array}
\right)_{i_1=1,\dots,m,~i_2=1,\dots,n}.
\]
\end{rem}
In the following we assume that the Runge-Kutta method is A-stable and $\Re(\varkappa_{kl})\leq0$
 or $|\kappa_{kl}|\leq C_2$ for all $h\in(0,h_0]$ with a positive constant $C_2$. Then for sufficiently small $\tau$
the matrix $G(\tau R_{kl})$ is regular, and the Runge-Kutta system \Nr{NLRKGlgWK1} has a unique solution.
 Using (\ref{NLRKGlg1WK1}), (\ref{Stern1}) and the transformed components $\delta_{1kl,m+1}$ of (\ref{ResGlg2})
  we obtain  the recursion
\begin{alignat}1
e_{1kl,m+1}
&=R(\tau R_{kl})e_{1kl,m}
+L(\tau R_{kl})\Delta_{1kl,m+1}
+\tau J(\tau R_{kl})\alpha_{h1kl, m+1}
+\delta_{1kl,m+1}\label{Fehlerrekursion}
\end{alignat}
for the discretization
 error $e_{1kl,m+1}$, where we have used the abbreviations
\[
J(z)=b^{\top}G(z)^{-1},
\quad
R(z)=1+J(z){1\kern-0.25em{\rm l}}_sz,
\quad
L(z)=J(z)z
\]
($R(z)$ equals the classical stability function
of the Runge-Kutta method).

Solving the recursion (\ref{Fehlerrekursion}) with $e_0=0$ leads to
\begin{alignat}1\nonumber
e_{1kl,m+1}
=&
\sum_{i=0}^mR(\tau R_{kl})^iL(\tau R_{kl})\Delta_{1kl,m+1-i}
\\&\label{emp1glgn1}
+\tau\sum_{i=0}^mR(\tau R_{kl})^iJ(\tau R_{kl})\alpha_{h1kl, m+1-i}+\sum_{i=0}^mR(\tau R_{kl})^i\delta_{1kl,m+1-i}.
\end{alignat}
Now we assume that the Runge-Kutta method under consideration has
(classical) order  $p$ and
 stage order $q$ ($p\geq q$). Then
 the simplifying conditions (see Hairer/Wanner \cite{HairerWanner})
 \begin{alignat*}1
B(p):~&\sum_{i=1}^sb_ic_i^{k-1}=\frac{1}{k},~k=1,\dots,p,\\
C(q):~&\sum_{j=1}^sa_{ij}c_j^{k-1}=\frac{1}{k}c_i^{k},~i=1,\dots,s,~k=1,\dots,q,
\end{alignat*}
are fulfilled.\\
With a Taylor expansion of $U_{h}(t_m+c_j\tau)$ and
 $\dot{U}_{h}(t_m+c_j\tau),~j=1,\dots,s$, around $t_m$ up to the order  $p$
we obtain for the $j$-th component of the residual error
$\Delta_{m+1}$ the equation
\begin{equation*}
\Delta_{j,m+1}
=
\sum_{r=q+1}^p
\frac{\tau^r}{r!}
\left(
\tilde{c}^r
-r\mathfrak{A}\tilde{c}^{r-1}
\right)_j
U_{h}^{(r)}(t_m)
+r_{\Delta_j,m+1},\quad\|r_{\Delta_j,m+1}\|=\mathcal{O}(\tau^{p+1})
\end{equation*}
with $\tilde{c}^i=(c_1^i,\dots,c_s^i)^\top$. Therefore, with
$r_{\Delta1kl,m+1}=\left(
r_{\Delta_11kl,m+1}^\top,\dots,r_{\Delta_s1kl,m+1}^\top\right)^\top$
and
\begin{equation*}
W_r(z)=\frac{L(z)\left[\tilde{c}^r-r\mathfrak{A}\tilde{c}^{r-1}\right]}{1-R(z)},
\end{equation*}
the error equation (\ref{emp1glgn1}) can be written as
\begin{alignat}1\nonumber
e_{1kl,m+1}~=~&
\tau\sum_{i=0}^mR(\tau R_{kl})^iJ(\tau R_{kl})\alpha_{h1kl, m+1-i}+\sum_{i=0}^mR(\tau R_{kl})^i\delta_{1kl,m+1-i}\\
 &
+\nonumber
\sum_{i=0}^mR(\tau R_{kl})^iL(\tau R_{kl})r_{\Delta1kl,m+1}\\
&
+\underbrace{\sum_{i=0}^mR(\tau R_{kl})^i
\left(
I_{n_{kl}}-R(\tau R_{kl})
\right)
\sum_{r=q+1}^p
\frac{\tau^{r}}{r!}
W_{r}(\tau R_{kl})U_{h1kl}^{(r)}(t_{m-i})}_{=\kappa}.\label{emglg2}
\end{alignat}
\begin{rem}
The function $W_{r}(z)$ was
introduced by Ostermann/Roche \cite{OstermannRoche} to investigate the convergence
of Runge-Kutta methods for abstract scalar parabolic differential equations.
\end{rem}
For the subsequent error estimate, the term $\kappa$ is transformed in the following manner: By exchanging the order
of summation we get
\begin{alignat*}1
\kappa
\nonumber
=&~\nonumber
\sum_{r=q+1}^p
\frac{\tau^{r}}{r!}
\Bigg(
\sum_{i=0}^{m-1}
R(\tau R_{kl})^{m-i}
W_{r}(\tau R_{kl})U_{h1kl}^{(r)}(t_{i})
+
W_{r}(\tau R_{kl})U_{h1kl}^{(r)}(t_{m})
\\&~\nonumber
-
\sum_{i=1}^{m}
R(\tau R_{kl})^{m-i+1}
W_{r}(\tau R_{kl})U_{h1kl}^{(r)}(t_{i})
-
R(\tau R_{kl})^{m+1}
W_{r}(\tau R_{kl})U_{h1kl}^{(r)}(t_{0})
\Bigg).
\end{alignat*}
From this we obtain
\begin{alignat*}1
\kappa
=&~\nonumber
\sum_{r=q+1}^p
\frac{\tau^{r}}{r!}
\Bigg(
\sum_{i=0}^{m-1}
R(\tau R_{kl})^{m-i}
W_{r}(\tau R_{kl})
\left( U_{h1kl}^{(r)}(t_{i})-U_{h1kl}^{(r)}(t_{i+1})\right)
\\&~\nonumber
+
W_{r}(\tau R_{kl})U_{h1kl}^{(r)}(t_{m})
-
R(\tau R_{kl})^{m+1}
W_{r}(\tau R_{kl})U_{h1kl}^{(r)}(t_{0})
\Bigg).
\end{alignat*}
Therefore it holds
\begin{alignat*}1
\kappa
=&~\nonumber
\sum_{r=q+1}^p
\frac{\tau^{r}}{r!}
\Bigg(
-\sum_{i=0}^{m-1}
R(\tau R_{kl})^{m-i}
W_{r}(\tau R_{kl})
\int\limits_{t_i}^{t_{i+1}}
U_{h1kl}^{(r+1)}(s)~ds
\\
&~ + W_{r}(\tau R_{kl})U_{h1kl}^{(r)}(t_{m}) - R(\tau
R_{kl})^{m+1} W_{r}(\tau R_{kl})U_{h1kl}^{(r)}(t_{0}) \Bigg).
\end{alignat*}
A similar transformation can be found in
Brenner/Crouzeix/Thom\'{e}e \cite{BrennerCrouzeixThomee}.

Inserting this into \Nr{emglg2} results in
\begin{alignat*}1\nonumber
e_{1kl,m+1}~=~&
\tau\sum_{i=0}^m
R(\tau R_{kl})^iJ(\tau R_{kl})\alpha_{h1kl, m+1-i}
+\sum_{i=0}^mR(\tau R_{kl})^i\delta_{1kl,m+1-i}\\
 &
+\nonumber
\sum_{i=0}^mR(\tau R_{kl})^iL(\tau R_{kl})r_{\Delta1kl,m+1}\\
&+
\nonumber
\sum_{r=q+1}^p
\frac{\tau^{r}}{r!}
\Bigg(
-\sum_{i=0}^{m-1}
R(\tau R_{kl})^{m-i}
W_{r}(\tau R_{kl})
\int\limits_{t_i}^{t_{i+1}}
U_{h1kl}^{(r+1)}(s)~ds
\\
&
+
W_{r}(\tau R_{kl})U_{h1kl}^{(r)}(t_{m})
-
R(\tau R_{kl})^{m+1}
W_{r}(\tau R_{kl})U_{h1kl}^{(r)}(t_{0})
\Bigg).
\end{alignat*}
Assuming that the Runge-Kutta matrix $\mbox{$\mathfrak{A}$}$ is re\-gu\-lar we can derive
 an analogous equation for the components $e_{2kl,m+1}$ of the transformed global discretization error
\begin{alignat*}1\nonumber
e_{2kl,m+1}~=~&
\tau\sum_{i=0}^m
\tilde{R}(N_{kl})^i\tilde{J}(N_{kl})\alpha_{h2kl, m+1-i}
+\sum_{i=0}^m\tilde{R}(N_{kl})^i\delta_{2kl,m+1-i}\\
 &
+\nonumber
\sum_{i=0}^m\tilde{R}(N_{kl})^i\tilde{L}(N_{kl})r_{\Delta2kl,m+1}\\
&+
\nonumber
\sum_{r=q+1}^p
\frac{\tau^{r}}{r!}
\Bigg(
-\sum_{i=0}^{m-1}
\tilde{R}(N_{kl})^{m-i}
\tilde{W}_{r}(N_{kl})
\int\limits_{t_i}^{t_{i+1}}
U_{h2kl}^{(r+1)}(s)~ds
\\
&
+
\tilde{W}_{r}(N_{kl})U_{h2kl}^{(r)}(t_{m})
-
\tilde{R}(N_{kl})^{m+1}
\tilde{W}_{r}(N_{kl})U_{h2kl}^{(r)}(t_{0})
\Bigg)
\end{alignat*}
with the abbreviations
\begin{alignat*}3
\tilde{J}(z)&=b^{\top}(I_sz-\tau\mbox{$\mathfrak{A}$})^{-1},\qquad&\tilde{R}(z)&=1+\tau\tilde{J}(z){1\kern-0.25em{\rm l}}_s,\\
\tilde{L}(z)&=\tau\tilde{J}(z),\qquad&
\tilde{W}_r(z)&=
\frac{\tilde{L}(z)\left[\tilde{c}^r-r\mathfrak{A}\tilde{c}^{r-1}\right]}{1-\tilde{R}(z)}.
\end{alignat*}
Finally, using (\ref{alpha2}), we get for
$e_{Qk,m+1}=Q_k(e_{1k1,m+1}^\top,\dots,e_{1ks_k,m+1}^\top,e_{2k1,m+1}^\top,\dots,e_{2kl_k,m+1}^\top)^\top$
the equation
\begin{alignat}1
e_{Qk,m+1}=&\nonumber
 h^{p_x}\sum_{j=1}^s
\tau \sum_{i=0}^m
Q_{k}
\mathrm{diag}\{
\dots,R(\tau R_{k j_1})^iJ_j(\tau R_{k j_1}),\dots,
\\&\nonumber
\hspace{3cm}\tilde{R}(N_{m_{k j_2}})^i\tilde{J}_j(N_{m_{k j_2}}),\dots
\}
P_{k}B\gamma_{hQk}(t_{m-i}+c_j\tau)\\\nonumber
&+\sum_{i=0}^m
Q_{k}
\mathrm{diag}\{
\dots,R(\tau R_{k j_1})^i,\dots,\tilde{R}(N_{m_{k j_2}})^i,\dots
\}
Q_{k}^{-1}
\delta_{Qk,m+1-i}
\\&\nonumber
+\sum_{j=1}^s
\sum_{i=0}^m
Q_{k}
\mathrm{diag}\{
\dots,R(\tau R_{k j_1})^iL_j(\tau R_{k j_1}),\dots,
\\&\nonumber
\hspace{3cm}\tilde{R}(N_{m_{k j_2}})^i\tilde{L}_j(N_{m_{k j_2}}),\dots
\}
Q_{k}^{-1}
r_{\Delta_jQ,m+1}\\\nonumber
&
+\sum_{r=q+1}^p
\frac{\tau^{r}}{r!}
Q_k\Bigg(
\mathrm{diag}_k\{\dots,
W_{r}(\tau R_{k j_1}),\dots,
\tilde{W}_{r}(N_{m_{k j_2}}),\dots
\}
Q_{k}^{-1}
U_{Qk}^{(r)}(t_m)
\\\nonumber&
-\sum_{i=0}^{m-1}
\int\limits_{t_i}^{t_{i+1}}
\mathrm{diag}_k\{\dots,
R(\tau R_{k j_1})^{m-i}
W_{r}(\tau R_{k j_1})
,\dots,
\tilde{R}(N_{m_{k j_2}})^{m-i}\tilde{W}_{r}(N_{m_{k j_2}}),\dots
\}\\\nonumber
&\hspace{3cm}
Q_{k}^{-1}
U_{Qk}^{(r+1)}(s)~ds
\\\nonumber
&
-\mathrm{diag}_k\{\dots,
R(\tau R_{k j_1})^{m+1}W_{r}(\tau R_{k j_1}),\dots,
\tilde{R}(N_{m_{k j_2}})^{m+1}\tilde{W}_{r}(N_{m_{k j_2}}),\dots
\}\\
&\hspace{3cm}
\label{emglg5}
Q_{k}^{-1}
U_{Qk}^{(r)}(t_0)
\Bigg).
\end{alignat}
Now we can estimate the different terms in \Nr{emglg5}. For that
purpose  we assume in the following that
 the matrix norms
\begin{subequations}\label{QkQkhm1QkPkBvBed}
\begin{equation}\label{QkQkhm1Bed}
\|Q_k\mathrm{diag}\{N^i_{n_{k1}},\mathfrak0,\dots,\mathfrak0\}Q_k^{-1}\|,\dots,\|Q_k\mathrm{diag}\{\mathfrak0,\dots,N^i_{m_{k l_k}}\}Q^{-1}_k\|
\end{equation}
and
\begin{equation}\label{QkPkBvBed}
\|Q_k\mathrm{diag}\{N^i_{n_{k1}},\mathfrak0,\dots,\mathfrak0\}P_k B\|,\dots,\|Q_k\mathrm{diag}\{\mathfrak0,\dots,\mathfrak0,N^i_{m_{k l_k}}\}P_k B\|
\end{equation}
\end{subequations}
are bounded for  $i=0,\dots,\max\{\nu_{dt},n_{kj_1}:~j_1=1,\dots,s_k\}-1$ and all $h\in(0,h_0]$, where $\mathfrak0$ denotes a zero matrix.\\
Because of the $A$-stability of the Runge-Kutta method and
$\Re(\varkappa_{k j_1})\leq0$ or $|\varkappa_{kj_1}|\leq C_2$
 for all $h\in(0,h_0]$ we have that
 $\|R(\tau R_{kj_1})^i\|$, $\|J_j(\tau R_{k j_1})\|$ and $\|L_j(\tau R_{k j_1})\|$ are bounded for sufficiently small $\tau$.
We assume further that $R(it)\neq1$ for $t\in\Rset\setminus\{0\}$ and $\lim\limits_{z\to-\infty}R(z)\neq1$. Then
$W_r(\tau R_{kj_1})$ exists and is bounded. Moreover, as it is shown in Ostermann/{\hspace{0mm}}Roche \cite{OstermannRoche},
 one has
  \[\|\tau^rW_r(\tau
  R_{kj_1})\|=\mathcal{O}(\tau^{\min\{p,q+2+\alpha\}})\|R_{kj_1}^{\max\{0,\min\{p-r,q+2+\alpha-r\}\}}\|\]
 with
 $\alpha\in\Rset$, $\alpha\geq-1$.

Assuming that $|\varkappa_{k j_1}|\leq C_3(1+|\lambda_k|)$ one can show (cf. D.\cite{Diss}, the proof
relies on the Mean Value Theorem and Abel's partial summation formula) that
\begin{equation*}
\|R_{kj_1}^{1+\alpha}\|~\|U_{Qk}^{(r)}(t_m)\|=k^{1+2\alpha}\mathcal{O}(h^{-\frac12}).
\end{equation*}
Altogether the terms in \Nr{emglg5} that originate from
$e_{1kl,m+1}$ are of order
\begin{equation}\label{GlgOrdNr1}
\mathcal{O}(h^{p_x})+\mathcal{O}(\tau^{\min\{p,q+2+\alpha\}})h^{-\frac12}k^{1+2\alpha}.
\end{equation}
With
the Taylor expansion
\[
\gamma_{hQk}(t_{m-i}+c_j\tau)
=\sum_{l=0}^{\nu_{dt}-2}\sum_{k=0}^{l}\frac{\tau^{l}(-i)^{l-k}}{k!(l-k)!}c_j^k\frac{\partial^l}{\partial t^{l}}\gamma_{hQk}(t_m)+\mathcal{O}(\tau^{\nu_{dt}-1}),
\]
the term
\begin{alignat*}1
a=&h^{p_x}\sum_{j=1}^s\tau
 \sum_{i=0}^m
Q_{k}
\mathrm{diag}\{0,
\dots,0,
\tilde{R}(N_{m_{k j_2}})^i\tilde{J}_j(N_{m_{k j_2}}),0,\dots,
\}
P_{k}B\gamma_{hQk}(t_{m-i}+c_j\tau)
\end{alignat*}
of $e_{Qk,m+1}$ can be written as
\begin{alignat*}1
a=&\sum_{j=0}^{\nu_{dt}-1}\frac1{j!}\sum_{l=0}^{\nu_{dt}-2}
Q_{k}
\mathrm{diag}\{0,
\dots,0,
\frac{N_{m_{k j_2}}^j}{\tau^{j-l}}\sum_{k=0}^l\frac1{k!(l-k)!}
\\&\nonumber
\qquad\qquad\qquad
\sum_{i=0}^m(-i)^{l-k}\left(\tau^{j+1}\tilde{R}(z)^i\tilde{J}(z)\tilde{c}^k\right)^{(j)}(0),
0,\dots,0
\}
P_{k}B\frac{\partial}{\partial t^{l}}\gamma_{hQk}(t_m)
\\&\nonumber
+\sum_{j=0}^{\nu_{dt}-1}\frac1{j!}
\sum_{i=0}^m
Q_{k}
\mathrm{diag}\{0,
\dots,0,
N_{m_{k j_2}}^j
\left(\tau^{j+1}\tilde{R}(z)^i\tilde{J}(z)\right)^{(j)}(0),0,\dots,0
\}
P_{k}B
\mathcal{O}(\tau^{\nu_{dt}-1})
,
\end{alignat*}
where $(\dots)^{(j)}$ denotes the $j$-th derivative w.r.t. $z$.\\
For $L$-stable Runge-Kutta  methods with regular coefficient matrix $\mbox{$\mathfrak{A}$}$ we have
 $R(\infty)=1-b^\top\mbox{$\mathfrak{A}$}^{-1}{1\kern-0.25em{\rm l}}=0$ and therefore
 $\tilde{R}(0)=0$.
 If the matrix norms in \Nr{QkQkhm1Bed}
are bounded,
 then $\|a\|$ is bounded if
 \begin{equation}\label{BedRhochiJalphabeschr}
\sum_{k=0}^l\frac1{k!(l-k)!}
\sum_{i=0}^j(-i)^{l-k}\left(\tau^{j+1}\tilde{R}(z)^i\tilde{J}(z)\tilde{c}^k\right)^{(j)}(0)=0
\end{equation}
for $l=0,\dots,j-1,~j=1,\dots,\nu_{dt}-1$.\\
The remaining terms in the equation \Nr{emglg5} yield
the classical order $p_{\nu_{dt}}$ of the Runge-Kutta method applied to a linear DAE of index $\nu_{dt}$
 with constant coefficients. Thus, altogether we have
 \[
e_{Qk,m+1}=\mathcal{O}(h^{p_x})+\mathcal{O}(\tau^{p_{\nu_{dt}}})+\mathcal{O}(\tau^{\min\{p,q+2+\alpha\}})h^{-\frac12}k^{1+2\alpha}.
 \]
From \Nr{eqzsglg} it follows that we have to choose $\alpha$ such that $\sum\limits_{k=1}^Nk^{2(1+2\alpha)}$ is
bounded for $N\to\infty$. This implies $\alpha=-\frac34-\varepsilon$, $\varepsilon>0$, and we have
\[
\|e_{m+1}\|=\mathcal{O}(h^{p_x})+\mathcal{O}(\tau^{\min\{p_{\nu_{dt}},q+1.25-\varepsilon\}}).
\]
If the derivatives of order $(q+1)$  w.r.t. the time of the
boundary conditions that enter into the space discretization
 are homogeneous, i.e.
\begin{equation}\label{qp1Ablhomogen}
B\frac{\partial^{q+1}u}{\partial t^{q+1}}=0
,
\end{equation}
(\ref{alphaWK1}) yields
\[
R_{kl}U_{h1kl}^{(r)}(t)={U}_{h1kl}^{(r+1)}(t)-\Tilde{F}_{1kl}(t)+\alpha_{h1kl}(t),
\]
and instead of (\ref{GlgOrdNr1}) we obtain the order
\[
\mathcal{O}(h^{p_x})+\mathcal{O}(\tau^{\min\{p,q+2+\alpha\}})h^{-\frac12}k^{2\alpha-1},
\]
which implies $\alpha=\frac14-\varepsilon$, and therefore
\[
\|e_{m+1}\|=\mathcal{O}(h^{p_x})+\mathcal{O}(\tau^{\min\{p_{\nu_{dt}},q+2.25-\varepsilon\}}).
\]
Summarized, we have the following convergence result for smooth
enough solutions $u(t,x)$ of the PDAE ($\max\{p+1,p+\nu_{dt}-1\}$
times differentiable with respect to $t$ in $[t_0,t_e]$ and
$p_x+2$ times differentiable with respect to $x$ in $[-l,l]$):
\begin{thm}\label{RKVWK}
Let the following assumptions be fulfilled for
 $h\to0$ \mbox{($N\to\infty$)} and $k=1,\dots,N$:
\begin{enumerate}\renewcommand{\labelenumi}{\alph{enumi})}
\item for the matrix pencils $D_k+\lambda A$ there exist Weierstrass-Kronecker decompositions
 according to \Nr{WKTransf}, and the matrix norms in \Nr{QkQkhm1QkPkBvBed} are bounded,
\item
$\Re(\varkappa_{k j_1})\leq0$ or $|\varkappa_{k j_1}|\leq C_2$ for all $h\in(0,h_0]$,
\item $|\varkappa_{k j_1}|\leq C_3(1+|\lambda_k|)$,
\item if $\nu_{dt}>2$ then \Nr{BedRhochiJalphabeschr} is  fulfilled for $l=0,\dots,j-1,~j=1,\dots,\nu_{dt}-1$.
\end{enumerate}
Furthermore let the Runge-Kutta method be of consistency order $p$, stage order $q$ and $L$-stable (if $\nu_{dt}=0$ or  $\nu_{dt}=1$, it suffices $A$-stability with
 $\lim\limits_{z\to-\infty}R(z)<1$) with a regular matrix $\mbox{$\mathfrak{A}$}$ and
  $R(it)\neq1$ for $t\in\Rset\setminus\{0\}$. Let $p_{\nu_{dt}}$ be the classical order of the Runge-Kutta method applied to a linear DAE of index $\nu_{dt}$
 with constant coefficients.\\
Then the discretization method (\ref{NLRKGlg})
 converges for linear PDAEs after $\nu_{dt}$ time steps
with the order $(p_x,p^\star)$ in the discrete $L_2$-norm in
space and in the maximum norm in time with
 \[
 p^\star=\left\{
\begin{array}{l@{:~}l}
\min\{p_{\nu_{dt}},q+1.25-\varepsilon\}&\text{inhomog. boundary conditions according to (\ref{qp1Ablhomogen})}\\
\min\{p_{\nu_{dt}},q+2.25-\varepsilon\}&\text{homog. boundary conditions according to (\ref{qp1Ablhomogen})}
\end{array}
\right.
\]
and $\varepsilon>0$ arbitrary small.
\end{thm}
\begin{rem}\label{BemKonvSatz}
\begin{enumerate}
\item Stage order  $q\geq\nu_{dt}-2$ implies condition d) in
Theorem 7 for $\nu_{dt}=3$ or $\nu_{dt}=4$.
\item The assumptions on the Runge-Kutta method are fulfilled, e.g., for the
Radau IIA and the Lobatto IIIC methods and in the case of $\nu_{dt}\leq1$ also for the implicit midpoint rule.
 \item If $p\geq p^\star+1$, then for $L$-stable Runge-Kutta methods with $zL(1,z)$ bounded for $\Re(z)\leq0$, the condition that
 \Nr{QkQkhm1Bed} is bounded can be replaced by the boundedness of the matrix norms
 \begin{subequations}\label{QkQkhm1Bed2}
 \begin{alignat}1
&\|\frac1{\kappa_{k1}}Q_k\mathrm{diag}\{N^i_{n_{k1}},\mathfrak0,\dots,\mathfrak0\}Q_k^{-1}\|,\dots,
\|\frac1{\kappa_{k s_k}}Q_k\mathrm{diag}\{\mathfrak0,\dots,\mathfrak0,N^i_{n_{k s_k}},\mathfrak0,\dots,\mathfrak0\}Q_k^{-1}\|,\\
&\|Q_k\mathrm{diag}\{\mathfrak0,\dots,\mathfrak0,N^i_{m_{k 1}},\mathfrak0,\dots,\mathfrak0\}Q^{-1}_k\|,
\|Q_k\mathrm{diag}\{\mathfrak0,\dots,N^i_{m_{k l_k}}\}Q^{-1}_k\|.
\end{alignat}
\end{subequations}
\item If we choose $\varepsilon=0$, then we get for $p^\star<p_{\nu_{dt}}$
 \[
 \|e_{m+1}\|=\mathcal{O}(h^{p_x})+
\mathcal{O}\left(\sqrt{|\ln h|}\tau^{p^\star}\right).
\]
\end{enumerate}
\end{rem}
\begin{rem}\label{BemEuler}
For a given Runge-Kutta method, Theorem \ref{RKVWK} can be specialized. E.g., if we take the implicit Euler method,
 the resulting BTCS method is convergent of time order 1 for arbitrary time index, if only the conditions a) and b)
 of Theorem \ref{RKVWK}  are fulfilled. For the Radau IIA methods with $s\geq2$ stages we get
 \[
 p^\star=\left\{
\begin{array}{l@{:~}l}
\min\{s+1.25-\varepsilon,2s-1\}&\nu_{dt}=0,1, \text{inhomog. b.c.s according to (\ref{qp1Ablhomogen})}\\
\min\{s+2.25-\varepsilon,2s-1\}&\nu_{dt}=0,1, \text{homog. b.c.s according to (\ref{qp1Ablhomogen})}\\
s+2-\nu_{dt}&\nu_{dt}\geq 2
\end{array}
\right.
\]
as temporal order of convergence, provided that the assumptions a)-d) of Theorem \ref{RKVWK} are fulfilled.
\end{rem}
\section{Numerical examples}\label{Abschnitt5}
The numerical examples given below illustrate our convergence
results. For the time integration we use the backward Euler
method and the  code RADAU5, which is  a variable step size
implementation of the 3-stage Radau IIA method, see Hairer/Wanner
 \cite{HairerWanner}. The Euler and Radau IIA methods are of great importance in
applications.
\begin{exmp}
The backward Euler method is given by the parameters
\[
s=1,\mbox{$\mathfrak{A}$}=(1),b=c=(1),p=q=1.
\]
We consider the linear PDAE
\[\underbrace{\begin{pmatrix}
       0&       1&       0\\
       0&       0&       1\\
       0&       0&       0\\
 \end{pmatrix}}_{=A}
 u_t+
 \underbrace{\left( \begin{array}{crr}
       0&       0&       -1\\
       0&       -1&       -1\\
       0&       0&       0\\
 \end{array} \right)}_{=B}
 u_{xx}+
\underbrace{ \left( \begin{array}{rrr}
       -1&       -1&       -1\\
       0&       -1&       0\\
       0&       0&       -1\\
 \end{array} \right)}_{=C}
 u=f,\]
 a coupled system of two parabolic equations and one algebraic
 equation, with
 $x\in [- 0.5,0.5]$, \mbox{$t \in [       0,       1]$}.
 The right-hand side, initial and Dirichlet boundary values are chosen such that
\[u(t,x)=
\left( \begin{array}{c}
 x(x-1)\sin(t)\\
 x(x-1)\cos(t)\\
 x(x-1)(e^t+t^5)\\
 \end{array} \right)\]
is the exact solution. It holds 
\[D_k=-{\lambda_k} B-C=\begin{pmatrix}
       1&       1~&      1+{\lambda_k}\\
       0&       1+{\lambda_k}&       {\lambda_k}\\
       0&       0~&       1
 \end{pmatrix}.\]
With \[P_k=\left( \begin{array}{crc}
       {\lambda_k}+1&       -1&       -{\lambda_k}-{\lambda_k}^2-1\\
       0&       1&       -{\lambda_k}-1\\
       0&       0&       1\\
 \end{array} \right)\text{ and }
Q_k=\begin{pmatrix}
       \frac1{{\lambda_k}+1}&       ~0&       ~0\\
       0&       \frac{~1}{{\lambda_k}+1}&       \frac{~1}{{\lambda_k}+1}\\
       0&       ~0&       ~1\\
 \end{pmatrix},
 \]
we obtain the Weierstrass-Kronecker decomposition
 \[
P_kAQ_k=\begin{pmatrix}
0&1&0\\0&0&1\\0&0&0
\end{pmatrix},~P_kD_kQ_k=\begin{pmatrix}
1&0&0\\0&1&0\\0&0&1
\end{pmatrix}.
 \]
Therefore, the PDAE has differential time index 3, and the assumptions \Nr{QkQkhm1QkPkBvBed} are fulfilled.
 Remark \ref{BemEuler} yields that the BTCS method is convergent after three steps of time order $1$.
 This is confirmed by the numerical experiment, Table \ref{Tabelle1} shows the observed order of convergence in
 time at $(x=1,t_e=1)$.
\begin{table}[ht]
\begin{center}
\begin{tabular}{c|cccccc}
$  0.1\tau^{-1}$&$2^{ 2}$ &$2^{ 3}$ &$2^{ 4}$ &$2^{ 5}$ &$2^{ 6}$
&$2^{ 7}$ \\ \hline $0.1h^{-1}$
\\$2^{   2}$&      0.81&      0.91&      0.96&      0.98&      0.99&      0.99
\\$2^{   3}$&      0.81&      0.91&      0.96&      0.98&      0.99&      0.99
\\$2^{   4}$&      0.81&      0.91&      0.96&      0.98&      0.99&      0.99
 \end{tabular}
 \caption{Numerically observed order of convergence  in the
discrete $L_2$-norm.\label{Tabelle1}}
\end{center}
 \end{table}
The notation of the first element 0.81 denotes the observed order
when refining the grid from $ (h=0.1/2^2,\tau=0.1/2)$ to
$(h=0.1/2^2,\tau=0.1/2^2)$, i.e., \mbox{$0.81=\log_{2} \xi$,}
where $\xi$ denotes the ratio of the error with
$(h=0.1/2^2,\tau=0.1/2)$ to the error with
$(h=0.1/2^2,\tau=0.1/2^2)$.  Furthermore, we see
that a simultaneous refinement of $h$ and $\tau$ yields no order
reduction.
\end{exmp}
\begin{exmp}\label{Bsp_KonvOrd4.25}
We consider the 3-stage Radau IIA method with
consistency order $p=5$ and stage order $q=3$, and the linear PDAE
\[
\underbrace{\left( \begin{array}{rrr}
0&2&0\\1&-1&0\\1&-1&0
 \end{array} \right)}_{=A}u_t+
\underbrace{\left( \begin{array}{rrr}
-1&0&0\\0&0&0\\0&0&-1
 \end{array} \right)}_{=B}u_{xx}+
\underbrace{\left( \begin{array}{rrr}
0&0&0\\0&-1&0\\0&0&1
 \end{array} \right)}_{=C}u
 =
 f(t,x)
\]
with
$x\in[-1,1]$ and $t\in[0,1]$.
 This example shows the dependence of the time order on the boundary values.
 \begin{enumerate}\renewcommand{\labelenumi}{\arabic{enumi}.}
 \item
We choose the right-hand side such that
\[u(t,x)=\left(
x^2e^{-t},x^2e^{-\frac12t},x^2\sin t
\right)
^\top
\]
 is the exact solution. Then we have inhomogeneous boundary values
\[u(t,\mp1)=\left(
e^{-t},e^{-\frac12t},\sin t
\right)
^\top
.\]
Furthermore it holds
\[D_k=\left(\begin{array}{crc}{\lambda_k}&0&0\\0&1&0\\0&0&{\lambda_k}-1\end{array}\right),~
P_kAQ_k=\begin{pmatrix}1&0&0\\0&1&0\\0&0&0\end{pmatrix},
~P_kD_kQ_k=\begin{pmatrix}\frac{2{\lambda_k}}{-{\lambda_k}-\eta_k}&0&0\\0&\frac{2{\lambda_k}}{-{\lambda_k}+\eta_k}&0\\0&0&1\end{pmatrix}\]
with
\[
P_k=
\begin{pmatrix}
\frac{{\lambda_k}-\eta_k}{4{\lambda_k}}&1&0\\
\frac{{\lambda_k}+\eta_k}{4{\lambda_k}}&1&0\\
0&\frac1{1-{\lambda_k}}&\frac1{{\lambda_k}-1}
\end{pmatrix},~
Q_k=
\begin{pmatrix}
\frac{4{\lambda_k}}{({\lambda_k}+\eta_k)\eta_k}&-\frac{4{\lambda_k}}{({\lambda_k}-\eta_k)\eta_k}&0\\
-\frac{{\lambda_k}}{\eta_k}&\frac{{\lambda_k}}{\eta_k}&0\\
-\frac{{\lambda_k}}{\eta_k({\lambda_k}-1)}&\frac{{\lambda_k}}{\eta_k({\lambda_k}-1)}&1
\end{pmatrix},~
\eta_k=\sqrt{{\lambda_k}^2+8{\lambda_k}}
\]
($|\lambda_k+8|>\frac12$ for $N>3$, i.e. $h<\frac12$).\\
The PDAE has therefore differential time index 1, and the conditions a)-c) of Theorem \ref{RKVWK} are fulfilled which
 yields convergence of time order $4.25-\varepsilon$.
  This is confirmed by the numerical experiment, see Table \ref{Tabelle2}.
\begin{table}[ht]
\centering
\begin{tabular}{c|cccc}
$  0.1\tau^{-1}$&$2^{ 1}$ &$2^{ 2}$ &$2^{ 3}$
 \\ \hline
$0.2h^{-1}$
\\$2^{   3}$&      4.27&      4.28&      4.30
\\$2^{   4}$&      4.26&      4.26&      4.26
\\$2^{   5}$&      4.26&      4.26&      4.25
\\$2^{   6}$&      4.26&      4.26&      4.25
 \end{tabular}
\caption{Numerically observed order of convergence  in the
discrete $L_2$-norm for inhomogeneous boundary values.\label{Tabelle2}}
\end{table}
\item
If instead the right-hand side is chosen such that
\[u(t,x)=
\left(
(x^2-1)e^{-t}+t^3\cos^2x,x^2e^{-\frac12t},(x^2-1)\sin t+t^3\sin^2x
\right)^\top
\]
 is the exact solution, then we have inhomogeneous boundary values
\[u(t,\mp1)=
\left(
t^3\cos^21,e^{-\frac12t},t^3\sin^21
\right)^\top
,\]
 but the derivatives of $Bu(t,x)$ of order 4 w.r.t. the time vanish. Therefore, we obtain
  the convergence order
 5 in time,
   see Table \ref{Tabelle4}.
 \begin{table}[ht]
\centering \begin{tabular}{c|ccc}
$  0.1\tau^{-1}$&$2^{ 2}$ &$2^{ 3}$ \\ \hline
$0.2h^{-1}$
\\$2^{   3}$&      5.00&      5.00
\\$2^{   4}$&      5.00&      5.00
\\$2^{   5}$&      5.00&      5.00
\\$2^{   6}$&      5.00&      5.00
 \end{tabular}
\caption{Numerically observed order of convergence  in the
discrete $L_2$-norm for inhomogeneous boundary values where the derivatives of $Bu$ of order 4 w.r.t.
 the time vanish.\label{Tabelle4}}
 \end{table}
\end{enumerate}
\end{exmp}
\begin{exmp}
We consider the 3-stage Radau IIA method
 and the linear PDAE \Nr{GlgSystemSupraMagnetSpule} describing the superconducting coil. It holds
 \[
D_k=\begin{pmatrix}
{\lambda_k}&-1&0&0&\\0&~\;{\lambda_k}&0&0&\\0&~\;0&1&0\\0&~\;0&0&1
\end{pmatrix},~P_kAQ_k=
\begin{pmatrix}
1&0&0&0\\0&1&0&0\\0&0&0&1\\0&0&0&0
\end{pmatrix},~P_kD_kQ_k=
\begin{pmatrix}
-\frac{i{\lambda_k}}{\sqrt{1-{\lambda_k}}}&0&0&0\\~\;0&\frac{i{\lambda_k}}{\sqrt{1-{\lambda_k}}}&0&0\\~\;0&0&1&0\\~\;0&0&0&1
\end{pmatrix}
 \]
 with
 \[
P_k=\begin{pmatrix}
-\frac{i}{\sqrt{1-{\lambda_k}}}&-\frac{i\sqrt{1-{\lambda_k}}}{{\lambda_k}}&1&-1\\
~\;\frac{i}{\sqrt{1-{\lambda_k}}}&~\;\frac{i\sqrt{1-{\lambda_k}}}{{\lambda_k}}&1&-1\\
~\;0&~\;0&1&-\frac1{\lambda_k}\\
~\;\frac1{\lambda_k}&~\;0&0&~\;0&
\end{pmatrix},~
Q_k=\begin{pmatrix}
\frac1{2(1-{\lambda_k})}&\frac1{2(1-{\lambda_k})}&\;~0&-\frac{\lambda_k}{1-{\lambda_k}}\\
\frac{\lambda_k}{2(1-{\lambda_k})}&\frac{\lambda_k}{2(1-{\lambda_k})}&\;~0&-\frac{\lambda_k}{1-{\lambda_k}}\\
-\frac{i{\lambda_k}}{2(1-{\lambda_k})^{\frac32}}&\;~\frac{i{\lambda_k}}{2(1-{\lambda_k})^{\frac32}}&-\frac{\lambda_k}{1-{\lambda_k}}&\;~0\\
-\frac{i{\lambda_k}^2}{2(1-{\lambda_k})^{\frac32}}&\;~\frac{i{\lambda_k}^2}{2(1-{\lambda_k})^{\frac32}}&-\frac{\lambda_k}{1-{\lambda_k}}&\;~0
\end{pmatrix}.
 \]
 The coil PDAE has therefore differential time index 2, and the conditions of Theorem \ref{RKVWK}
 (with the matrix norms \Nr{QkQkhm1Bed} replaced by \Nr{QkQkhm1Bed2})
  are fulfilled which yields an order of convergence in time of $3$.

 This is confirmed by the numerical experiment, see Table \ref{Tabelle6}.
 \begin{table}[ht]
\centering \begin{tabular}{c|ccc}
$  0.1\tau^{-1}$&$2^{ 4}$ &$2^{ 5}$ &$2^{ 6}$ \\ \hline
$0.2h^{-1}$
\\$2^{   2}$&      3.00&      3.00&      3.00
\\$2^{   3}$&      3.00&      3.00&      3.00
\\$2^{   4}$&      3.00&      3.00&      3.00
 \end{tabular}\caption{Numerically observed order of convergence  in the
discrete $L_2$-norm for the coil PDAE.\label{Tabelle6}}
\end{table}
 \end{exmp}\vspace{-0.048cm}
\section{Conclusion}\label{Abschnitt6}
The attention has here been restricted to a class of linear
partial differential-algebraic equations. We have given
convergence results in dependence on the type of boundary values
and the time index. When the error is measured in the
discrete $L_2$-norm over the whole domain, the convergence order in time of the Runge-Kutta method for a smooth
solution is in general non-integer and smaller than the
order  expected for differential-algebraic equations of the same index. Some numerical
examples were presented and
confirm the theoretical convergence results.\\
The extension of the analysis to the case of  space $d$ dimensional
linear partial differential-algebraic equations of the form
\begin{equation*} A\:u_t(t,\vec{x}) +\sum_{i=1}^d B_i\:\left(
u_{x_ix_i}(t,\vec{x})+r_iu_{x_i}(t,\vec{x})\right) + C\: u(t,\vec{x}) =
f(t,\vec{x})
\end{equation*}
with $\vec x=(x_1,\cdots,x_d)^\top$ and a cuboid as domain is
possible, see D. \cite{Diss}, but becomes rather technical and  offers no new insight.
Furthermore, the consideration of periodic boundary values is
also possible. Here we could show as temporal convergence order
the order of an ordinary differential-algebraic equation. In the
case of Neumann boundary conditions, the temporal convergence
order lies in between
the order obtained for Dirichlet- and the order obtained for periodic boundary conditions.\\
Future work in this area will be concerned with  convergence
investigations for semi-linear partial differential-algebraic
equations.

{\bf Acknowledgements}

The authors are very grateful to the referee for his comments and
fruitful suggestions.

\end{document}